\input amstex
\input psfig.sty

\documentstyle{amsppt}

\NoRunningHeads
\magnification\magstep 1
\baselineskip 15pt
\pagewidth{6.4 truein}
\pageheight{8.6 truein}

\NoBlackBoxes

\TagsOnRight

\catcode`\@=11
\redefine\logo@{}
\catcode`\@=13

\def\qed{\hfill\ \vbox{\hrule\hbox{\vrule\kern4pt\vbox{\kern4pt{}
\kern4pt}\kern4pt\vrule}\hrule}}

\hfill 10/12/99%

\topmatter
\title 
Mathematical Problems in the Control of Underactuated Systems
\endtitle

\author  {\ }\endauthor

\affil
{\smc David Auckly}\\
Department of Mathematics, \\
Kansas State University, Manhattan, KS 66506-2602, USA\\
email: dav\@math.ksu.edu\\
\\
{\smc Lev Kapitanski} \\
Department of Mathematics,   \\
Kansas State University, Manhattan, KS 66506-2602, USA\\
email: levkapit\@math.ksu.edu \\
\endaffil

\thanks This work was partially supported 
by a grant from the National Science Foundation 
CMS 9813182.
\endthanks

\endtopmatter


There are many interesting mathematical problems in control theory. 
In this paper we will discuss problems and techniques related to 
underactuated systems. An underactuated system is one with fewer 
control inputs than degrees of freedom. Balancing a ruler on the 
tip of a finger is a good example of an underactuated system. 
This system has five degrees of freedom (three for the fingertip and two 
angles for the ruler). However, only the three degrees of freedom 
for the fingertip are directly controlled. In fact, 
any system requiring balance is an underactuated system. 
A bicycle is an obvious example. An airplane is a less obvious example 
(six degrees of freedom, underactuated by two). 

We will give a mathematical formulation of several problems arising from applications, review some standard and new techniques, and pose 
some interesting and challenging open questions. 

\head Stabilization of underactuated systems \endhead

To describe a mechanical system we start with a manifold, $Q$, 
representing all possible configurations of the system.
The configuration space $Q$ is equipped with a Riemannian metric, 
$g$, so that the kinetic energy is $\frac 12 g(\dot x, \dot x)$. 
It also comes with a function $V:\,Q\to \Bbb R$, 
and two fiber preserving maps $c, f:\,TQ\to TQ$. 
The function $V$ represents potential energy, $c$ represents 
dissipation, and $f$ represents applied external forces. 
The equations of motion are given by 
$$ 
\nabla_{\dot\gamma}\dot\gamma 
+ c(\dot\gamma) +\ grad_\gamma V
   =f(\dot\gamma).  \tag1 
$$
The external forces, $f$, are used to control the system. 
The system obtained by setting the control input to zero 
is called {\it the open loop system}. Equation (1) including 
the control input is called {\it the closed loop system}.
The system is {\it underactuated} if $f$ is restricted to be $0$ 
in some directions. In other words, there is a $g$-orthogonal projection 
$P$ onto the subspace of unactuated directions, and $P(f)$ must 
vanish. 

The basic problem is to find a function $f$ in some class so that 
solutions to equation (1) have some desired properties. 
Physically one may not always be able to measure the full state 
$(x,\dot x)\in TQ$ of the mechanical system. 
In this situation the function $f$ must only 
depend on the observable variables. For now, we will consider the case
when all variables may be observed. 
This is referred to as full state feedback
control. 

The stabilization problem is to find a control input, so that
some point $(x_0,0)$ will be an asymptotically stable equilibrium. Other
notions of stability may also be considered, however, asymptotic stability
is the most useful in applications. 
We will next review several approaches to the stabilization problem.

The most commonly employed technique used to address this problem is
linearization. Choosing $f$ so that the eigenvalues of the linearized
equation lie in the left half
plane, will ensure that the desired point is a locally asymptotically
stable equilibrium. This reduces the problem to an algebraic question
that can be easily solved and implemented. 

First, note that it is not possible to stabilize every system, for example,
$$
\left\{
\aligned
\ddot x^1\,-\,x^1 &\,=\,0\\
\ddot x^2\,+\,x^2 &\,=\,u\,.
\endaligned
\right.
$$
For a linear $n\times n$ 
system of the form $\dot y = A y + B u$ necessary and sufficient 
conditions for the existence of a linear stabilizing control law 
$u=Cy$ 
are well known \cite{14, 18}: the rank of the matrix $[sI-A, B]$ must be 
$n$ for all $\text{Re}\,s\ge 0$. 
If a linear stabilizing control law exists, 
there is a finite dimensional family of linear stabilizing control laws. 

Once it is known that a  stabilizing control law exists, one must 
choose a specific control law. 
The problem of finding a matrix $C$ given the eigenvalues of $A+BC$ 
is called pole placement. Engineers use various rules of thumb to 
decide where to place the poles. These rules of thumb are based upon 
the behavior of solutions to a constant coefficient second order ODE. 
For higher order systems one purposefully places two 
dominant poles, $z_1$ and $z_2$, with the remaining poles near 
the real axis and far to the left. This enables one to approximate 
solutions of the higher order system by solutions of a second order 
system. 

\head Example 1: The inverted pendulum cart \endhead 
\bigskip

\hskip 100bp\psfig{file=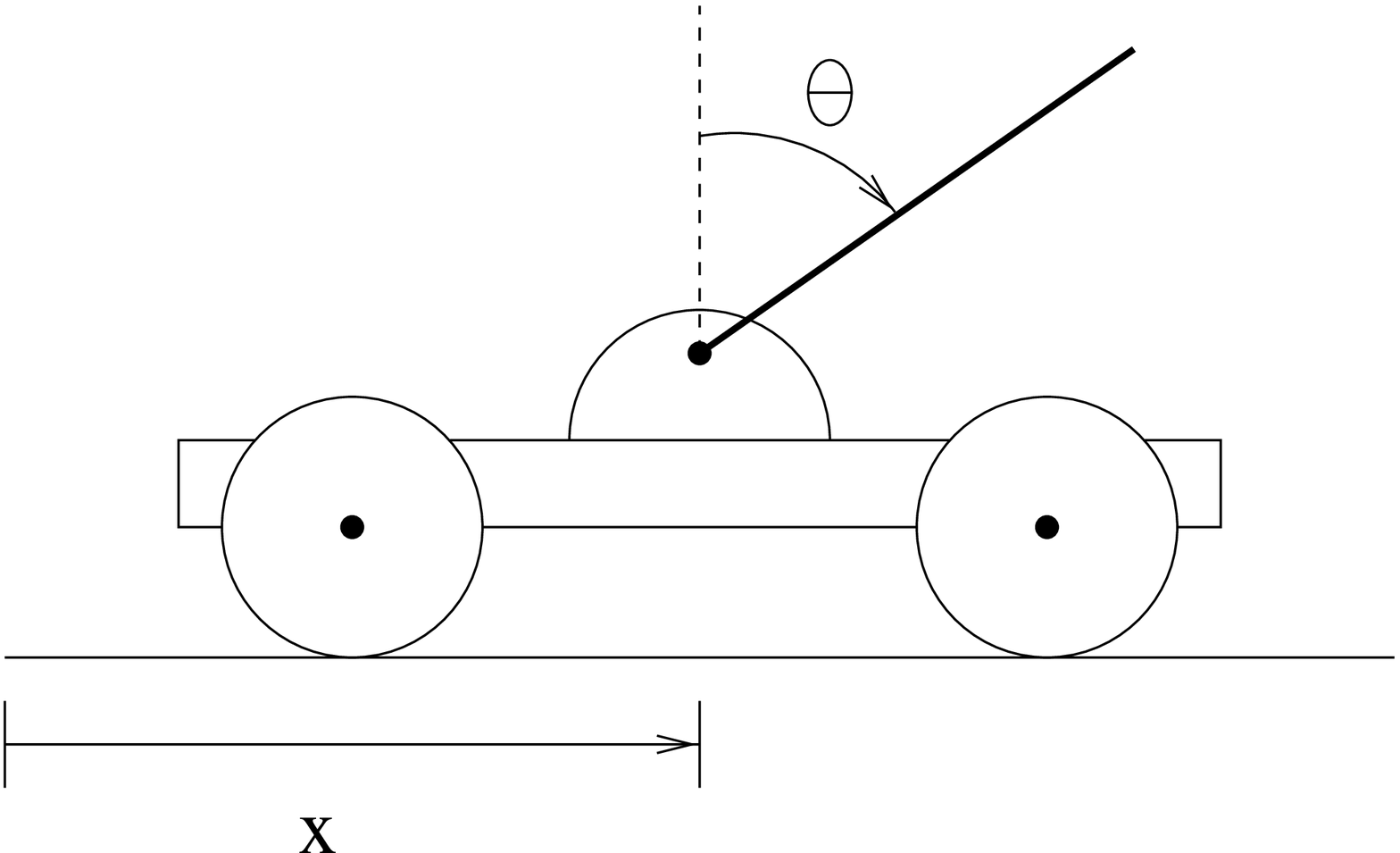,width=4truein,height=2.3truein}%

\medskip
\centerline{Figure 1}
\smallskip

With appropriate scaling the metric $g$ is given by 
$g=d\theta^2 +2 b \cos(\theta)\,dx\,d\theta + dx^2$, where $b$ is a physical parameter, $0<b<1$. The potential energy is given by $V=\cos(\theta)$. 
Since no torques can be applied directly to the pendulum, 
$P=(b\cos(\theta)\,dx+d\theta)\otimes {\partial/\partial\theta}$ is 
the orthogonal projection onto the direction ${\partial/\partial\theta}$. 
Assuming that there is no dissipation, $c=0$. 
The ${\partial/\partial\theta}$- and 
${\partial/\partial x}$-components of the equations of motion read: 
$$
\aligned 
\ddot\theta + b\cos(\theta)\,\ddot x -\sin(\theta) &= 0 \\    
b\cos(\theta)\,\ddot\theta + \ddot x -b\sin(\theta)\,\dot\theta^2 &= u, 
\endaligned
$$
where $u=g({\partial/\partial x},f)$ represents the external force 
applied to the base of the cart, and it is the control input. 

The linearization around $\theta=0, \,x=0,\,\dot\theta=0,\,\dot x=0$ 
reads
$$
{d\hfil\over dt}\,
\pmatrix \theta\\ x\\ \dot\theta\\ \dot x\endpmatrix\,=\,
\pmatrix 0 & 0 & 1 & 0\\
 0 & 0 & 0 & 1\\ {1-a_1 b\over 1-b^2} & {-a_2 b\over 1-b^2}  
& {-a_3 b\over 1-b^2} & {-a_4 b\over 1-b^2}\\
{a_{1}+b\over 1-b^2} & {a_{2}\over 1-b^2}  & {a_3\over 1-b^2} 
& {a_4\over 1-b^2}
\endpmatrix
\,
\pmatrix \theta\\ x\\ \dot\theta\\ \dot x\endpmatrix 
$$
where $\;u=a_1\theta+a_2 x+a_3 \dot\theta+a_4 \dot x\;$ is the 
linearized control input.
The characteristic polynomial of the above matrix is 
$$
\lambda^{4}+{a_3 b-a_4\over 1-b^2}\lambda^{3}
+{1+a_2-a_1 b\over 1-b^2}\lambda^{2}
+{a_4 (1+b^2)\over (1-b^2)^2}\lambda+{a_2 (1+b^2)\over (1-b^2)^2}.
$$
Once the desired eigenvalues are specified, it is an easy matter 
to solve for the $a_k$ and get the control law. 

\bigskip
Another standard technique employed in control design is linear quadratic 
optimal control, \cite{18}. For linearized systems $\dot y = Ay + Bu$, one looks for 
control laws which will minimize the functional 
$$
J(y,u)
=\int_0^\infty \langle y(t),Qy(t)\rangle +\langle u(t),Nu(t)\rangle\,dt 
$$
subject to the condition $\dot y = Ay + Bu$. Here $Q$ and $N$ are positive definite quadratic matrices which are usually specified by 
engineering rules of thumb.

If a stabilizing control law exists, the resulting 
control input can be expressed 
as a linear function of the state $y$. 

\bigskip

Linearization works very well for many practical applications. 
Unfortunately, the limitations of linearization are seldomly discussed. 
Continuing in this tradition we will now address nonlinear 
methods without stating why. 

Many nonlinear methods employ the notion of a Lyapunov function. 
For a system of ODEs, $\dot x=f(x)$, a Lyapunov function is a function 
$F(x)$ which is bounded from below and decreases along the trajectories, 
\cite{4}. 
It can be defined globally or locally. If $x_0$ is a stationary solution of 
$\dot x=f(x)$, and a unique local minimizer of the Lyapunov function $F$, 
then $x_0$ is Lyapunov stable. If $F$ is nonconstant along nonstationary 
trajectories then $x_0$ will be locally asymptotically stable. 
If, in addition, $F$ is defined globally, 
then $x_0$ is a global asymptotically stable equilibrium. 

For the stabilization problem, $\dot x=f(x,u)$, 
one wishes to find a control input, 
$u$, as a function of the state, $x$, and a Lyapunov function, 
$F(x)$, so that the closed loop system admits 
$F(x)$ as a Lyapunov function. In the linear case, 
$\dot y = Ay + Bu$,  with $u=Cy$, one looks for a quadratic 
Lyapunov function $F(y)=\langle y, K y\rangle$. The  
function $F$ will be a Lyapunov function if and only if 
$$
D\,=\,(A+BC)^*\,K + K\,(A+BC)\tag 2
$$
is negative definite. Given any negative definite $D$ one may solve 
equation (2) for $K$ if and only if all the eigenvalues of 
$A+BC$ are in the left half-plane. The solution is unique: 
$$
K = - \int_0^\infty e^{t(A+BC)^*}\,D\, e^{t(A+BC)}\,dt.
$$
Note, that in the optimal control approach discussed previously, 
the value function
$$
F(x)=\inf_{u}\,J(x,u)
$$
is a Lyapunov function.  
We also remark that, 
if there is a locally asymptotically stabilizing control law, 
then there exists a local Lyapunov function, \cite{4}.

Several recent papers propose to find control inputs so that 
the closed-loop system (1) would have a natural 
candidate for a Lyapunov function, \cite{1-3, 5-9, 13, 17}. 
In \cite{3} we introduce the following approach to the stabilization 
problem for underactuated systems. 

We consider the control problem 
$$ 
\nabla_{\dot\gamma}\dot\gamma 
+ c(\dot\gamma) +\ grad_\gamma V
   =f(\dot\gamma),  \tag1 
$$
subject to the constraint $P(f)=0$. 
We wish to find $f$ and a Lyapunov function simultaneously. 
In fact, we will describe an infinite dimensional family 
of control inputs.

Our approach to this question is to find functions 
$\widehat g$, $\widehat c$, 
$\widehat V$ and $f$ so that
solutions to Equation (1) are automatically solutions to
$$ 
\widehat \nabla_{\dot\gamma} \dot\gamma 
+ \widehat c(\dot\gamma) +\ \widehat  {grad}_\gamma
   \widehat V=0. \tag3 
$$
This is the matching philosohy. The motivation for this philosophy
is that $$\widehat H(X)=\frac{1}{2} \widehat g(X,X) + \widehat V$$
is a natural candidate for a Lyapunov function because
$\frac{d}{dt} \widehat H(X) = -\widehat g(\widehat c(X),X)$. A state,
$X_0 \in TQ$ will be an asymptotically stable equilibrium if
$\widehat H(X) \geq 0$, and $ -\widehat g(\widehat c(X),X) \geq 0$ 
with equality
only at $X_0$.

Equations (1) and (3) clearly hold if and only if:
$$ f(X)\equiv \nabla_X X-\widehat \nabla_X X+\ grad_\gamma V
   -\ \widehat {grad}_\gamma \widehat V +c(X)-\widehat c(X),
   \tag4 
$$
for every vector field $\,X$.
The condition $P(f)=0$ then becomes a system of nonlinear partial
differential equations for $\widehat g$, $\widehat V$, and 
$\widehat c$.
Notice that constant multiples of $g$, $c$, and $V$ satisfy
$\,P(f)=0$ even when $\,P$ has full rank. Thus, one would expect many
solutions when $\,P$ does not have full rank. Separating
$P(f)=0$ into terms which are quadratic in the velocity, independent
of the velocity or
odd functions of the velocity gives:
$$ 
P(\nabla_X X-\widehat \nabla_X X)=0, \tag5.1 
$$
$$ 
P(grad_\gamma V-\ \widehat  {grad}_\gamma \widehat V)=0, \tag5.2 
$$
$$ 
P(c(X)-\widehat c(X))=0. \tag5.3 
$$
We will look for solutions to these matching equations with
$\widehat g$ non-degenerate so that $g(X,Y)=\widehat g(\lambda X,Y)$ with
$\lambda\in\Gamma(T^\ast Q\otimes TQ)$. It is clear that $\lambda$
has to be
$g$ self-adjoint, i.e., $g(\lambda X,Y)=g(X,\lambda Y)$. 
We will derive a linear system of partial
differential equations for $\lambda$ which must be
satisfied if $\,\widehat g$ is to solve Equation (5.1).

In our previous paper, we described a method to find every solution to the 
matching equations by solving three linear systems  of partial differential
equations in a row. We will review this method now.

One first solves the equations 
$$ 
\nabla g\lambda\big|_{\text{Im}\ P^{\otimes 2}}=0,
   \tag6 
$$
for $\lambda |_{\text{Im}\ P}$. Then one solves 
$$
L_{{}_{\lambda PX}}\widehat g = L_{{}_{PX}} g \tag 7
$$
(this is a slight rewrite of equation (1.12) of our previous paper \cite{3}),
$$
L_{{}_{\lambda PX}}\widehat V = L_{{}_{PX}} V \tag 8
$$
(this is equation (1.13) of our previous paper \cite{3}),
then after solving equation (5.3), the control input will be given by (4). 

In the previous paper we explicitly showed that any solution 
to the matching equations solves equations (6), (7), (8), and (5.3) 
(Propositions 1.1, 1.2, and 1.3 of \cite{3}).  
Implicit in \cite{3, Proposition 1.4} is the fact that any solution 
of equations (6), (7), (8) and (5.3) is in turn a solution 
to the matching equation. We will make this argument explicit 
now.  Indeed, taking into account the fact that $P$ is $g$-selfadjoint, 
our $\widehat g$-Equation (7) implies the matching equation (5.1). 
Here is a short proof. For any $Z$ and $X$ we have
$$
\aligned
& g(P(\widehat\nabla_ZZ -\nabla_ZZ),\;X) \\
=& g(\widehat\nabla_ZZ -\nabla_ZZ,\;PX)\\
=& \widehat g(\widehat\nabla_ZZ,\;\lambda PX) - g(\nabla_ZZ,\;PX)\\
=& Z \widehat g(Z,\;\lambda PX) -\widehat g(Z,\;\widehat\nabla_Z \lambda PX)
- Z g(Z,\;PX) + g(Z,\;\nabla_Z PX)\\
=& - \widehat g(Z,\;\widehat\nabla_{\lambda PX} Z)+
     \widehat g(Z,\;[\lambda PX,Z])
+ g(Z,\;\nabla_{PX} Z) - g(Z,\;[PX,Z])\\
=& -\frac 12\, (L_{{}_{\lambda PX}}\,\widehat g(Z,Z))
+\frac 12\,(L_{{}_{PX}}\, g(Z,Z))\\
=& 0.
\endaligned
$$
Since this is true for all $X$, $\;P(\widehat\nabla_ZZ -\nabla_ZZ)=0$.
\medskip

Equations (6) and (7) imply additional compatibility conditions. 
Even though we do not know all the compatibility conditions in general, 
we do know all the compatibility conditions for systems with two 
degrees of freedom. Let us summarize our method in the case of 
two degrees of freedom one of which is unactuated.
Since the unactuated subspace is one dimensional, it can be locally 
expressed as the span of a unit length vectorfield, $PX$.  
Choose coordinates $x^1$, $x^2$ so that $PX={\partial\over\partial x^1}$. 
In these coordinates $g_{11}=1$. We will always write 
$\lambda PX=\sigma {\partial\over\partial x^1}
+\mu {\partial\over\partial x^2}$, where $\sigma$ and 
$\mu$ are yet to be found. 
The $\lambda$-equation may be rewritten as
$$
{\partial\over\partial x^1} (g_{11}\sigma+g_{12}\mu)-2[11,2]\,\mu=0,
\qquad
{\partial\over\partial x^2} (g_{11}\sigma+g_{12}\mu)-2[12,2]\,\mu=0. \tag9
$$
Here, 
$$[ij,k]=g(\nabla_{\partial_i}\partial_j,\partial_k)=
\frac12 (g_{ik,j}+g_{jk,i}-g_{ij,k}). \tag10
$$ 
For these equations to be consistent 
the following compatibility condition must hold:
$$
{\partial\over\partial x^2}([11,2]\,\mu)=
{\partial\over\partial x^1}([12,2]\,\mu). \tag11
$$
Notice that this is a first order partial differential equation 
equation. Theoretically, it can be solved for $\mu$ via the method of 
characteristics. Generically, a solution will include 
an arbitrary function of a single variable. 
Once $\mu$ is known, $\sigma$ is given by
$$
\sigma(x^1,x^2)=g_{12}(x^1,x^2)\mu(x^1,x^2)+
2\,\int \left([11,2]\,\mu(x^1,x^2)\,dx^1+
[12,2]\,\mu(x^1,x^2)\,dx^2\right)\,.
$$

The next step is to solve equations (6) for $\widehat g$. 
It turns out that it is easiest to solve first for $\widehat g_{11}$ 
and then find the remaining components from the algebraic system
$$
g=\widehat g \lambda.
$$
First note that the $\{11\}$ component of the right side of (7) is 
$$
(L_{PX}g)({\partial\over\partial x^1},{\partial\over\partial x^1})=
{\partial\over\partial x^1} g_{11}=0.
$$
Next, 
$$
(L_{\lambda PX}\widehat g)
({\partial\over\partial x^1},{\partial\over\partial x^1}) = 
{\lambda PX}(\widehat g_{11})
-2\widehat g([{\lambda PX},{\partial\over\partial x^1}],
{\partial\over\partial x^1}).
$$
Since
$$
[{\lambda PX},{\partial\over\partial x^1}]=
-{\partial\sigma\over\partial x^1} {\partial\over\partial x^1}
-{\partial\mu\over\partial x^1} {\partial\over\partial x^2} =
(-{\partial\sigma\over\partial x^1} {\partial\over\partial x^1}
+{\sigma\over\mu}{\partial\mu\over\partial x^1}){\partial\over\partial x^1}
-{1\over\mu}{\partial\mu\over\partial x^1}\lambda PX
$$
and $\widehat g \lambda = g$, we obtain, 
$$
\lambda PX \widehat g_{11}-2(-{\partial\sigma\over\partial x^1} 
+{\sigma\over\mu}{\partial\mu\over\partial x^1})\widehat g_{11}
+2{1\over\mu}{\partial\mu\over\partial x^1}g_{11}=0.
$$
Thus, $\widehat g_{11}$ satisfies the following first order 
partial differential equation
$$
\sigma {\partial\widehat g_{11}\over\partial x^1}+
\mu {\partial\widehat g_{11}\over\partial x^2}
+2({\partial\sigma\over\partial x^1}
-{\sigma\over\mu}{\partial\mu\over\partial x^1})\widehat g_{11}
+2{{\partial\mu\over\partial x^1}\over\mu}=0.
$$
The general solution to this first order PDE has an arbitrary function 
of a single variable in it. It is, once again, possible to solve 
this equation by the method of characteristics. Denote by $y( x^1,x^2)$ 
any solution of the homogeneous equation
$$
\sigma {\partial y\over\partial x^1}+
\mu {\partial y\over\partial x^2}=0 \tag 12
$$
such that $\partial y/\partial x^2\neq 0$. Let $\bar\sigma$ and $\bar\mu$ 
be $\sigma$ and $\mu$ considered as functions of $ x^1$ and $y$, i.e.,
$$
\bar\sigma( x^1,y( x^1,x^2))=\sigma( x^1,x^2),\qquad
\bar\mu( x^1,y( x^1,x^2))=\mu( x^1,x^2). 
$$
Then the solution to equation (7) is given explicitly by
$$
\widehat g_{11}( x^1,x^2)={\mu^2\over\sigma^2}\left[
-2\,\int_0^{x^1}{\bar\sigma\over\bar\mu^3}
{\partial\bar\mu\over\partial x^1}\,d x^1\,\big|_{y=y( x^1,x^2)}
+h(y( x^1,x^2))\right], \tag 13
$$
where $h(y)$ is an arbitrary function of a single variable.
After $\widehat g_{11}$ 
is found, we have
$$
\widehat g_{12} ={1\over\mu} (g_{11}-\sigma \widehat g_{11}),
\qquad 
\widehat g_{22} ={1\over\mu} (g_{12}-\sigma \widehat g_{12}). \tag 14
$$
The equation for $\widehat V$ reads
$$
\sigma {\partial\widehat V\over\partial x^1}+
\mu {\partial\widehat V\over\partial x^2}={\partial V\over\partial x^1}.
$$
Again, this equation can be solved by the method of characteristics.  
Once a solution $y$ of (12) is known, $\widehat V$ is given by: 
$$
\widehat V( x^1,x^2) = -\int_0^{x^1} {\bar V_{x^1}\over\bar\sigma}
\,d x^1\,\big|_{y=y( x^1,x^2)}+w(y( x^1,x^2))\,, \tag 15
$$
where $\bar V_{x^1}(x^1,y(x^1,x^2))={\partial V(x^1,x^2)\over\partial x^1}$.  
Also,  
$$
\widehat c_1=c_1+g_{12}\,(\widehat c_2-c_2).
$$
It is convenient to write the non-zero control input as a sum: 
$$
u=g(f,{\partial\over\partial x^2})=
g(\nabla_X X-\widehat \nabla_X X+\ grad_\gamma V
   - \widehat {grad}_\gamma \widehat V 
+c(X)-\widehat c(X),{\partial\over\partial x^2})=u_g+u_V+u_c, \tag 16
$$
where
$$
u_g=g(\nabla_X X-\widehat \nabla_X X,{\partial\over\partial x^2}),\;
u_V=g(grad_\gamma V
   - \widehat {grad}_\gamma \widehat V, 
{\partial\over\partial x^2}),\;
u_c=g(c(X)-\widehat c(X),{\partial\over\partial x^2}),
$$
In coordinates, 
$$
u_g=([ij,2]-g_{k2}\widehat\Gamma_{ij}^k)\dot x^i\dot x^j, \tag 17
$$
where $\widehat\Gamma_{ij}^k$ are the Christoffel symbols 
corresponding to $\widehat g$,
$$
\widehat\Gamma_{ij}^k=\widehat g^{kp}\,\widehat{[ij,p]},
$$
$\widehat g^{kp}$ is the inverse matrix to $\widehat g_{lm}$, and $\widehat{[ij,p]}$ is defined as in (9) with all $g$'s 
replaced by $\widehat g$'s.
The next term is: 
$$
u_V=V_{x^2}-\widehat g^{ij}\,g_{j2}\,\widehat V_{x^i}. \tag 18
$$
Finally, 
$$
u_c=\det g \,(c_2-\widehat c_2). \tag 19
$$
\smallskip

From the explicit formulae (16)-(19) one sees that the first order 
germs of our control inputs contain every possible linear control 
law: 
\proclaim{ Theorem} If a system with 2 degrees of freedom 
underactuated by 1 is linearly stabilizable, then, by choosing 
appropriate solutions of our matching equations, the linearization 
of the controlled system will have prescribed eigenvalues 
in the left halfplane. Even if the system is not stabilizable, 
any linear control law may be obtained as the first order germ 
of some control law in our family.
\endproclaim
\medskip

\head {Open problems}\endhead
\medskip

One of the reasons that the limitations of linearization 
are seldomly discussed is that there is no effective 
criteria to compare arbitrary control laws. To be more specific 
we will restrict our discussion to the stabilization problem.

A good stabilizing control law will produce a large basin of 
attraction, send solutions to the equilibrium in a short period 
of time, and will have low cost. Of course, the precise meaning 
of ``large", ``short", and ``low" depends on the concrete 
engineering problem. It is not clear how to quantify 
these concepts. 
For linear systems the size of the basin of attraction
 is irrelevant since 
the whole space is the basin of attraction of a stable equilibrium. 
For nonlinear systems this question is subtle.
One could just use the volume or diameter as a  measure 
of the size of the basin of attraction. 
These are not, however, usually appropriate measures of size, 
see Figure 2. In addition, they are difficult to compute. 

\bigskip


\hskip 60bp\psfig{file=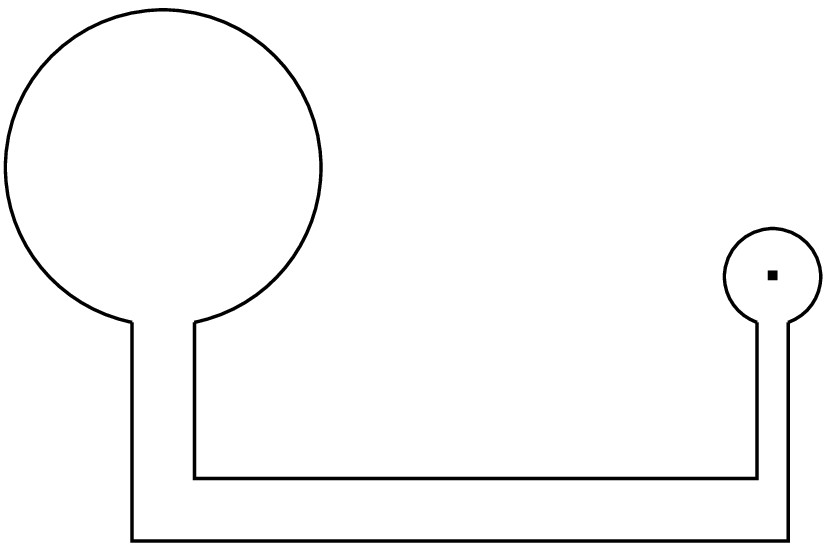,width=4truein,height=2.3truein}%

\centerline{Figure 2}
\smallskip

Alternatively, one could measure the radius of the largest 
inscribed ball centered at the equilibrium. It is not, however, clear 
which metric should be used in state space.
It is difficult to analytically or numerically estimate this radius. 
Even if one could compute this radius, it might not be the most 
relevant measure of performance. Real systems have an operating range:
rollerblades are not designed to handle Mach 2. 
So, let 
$\Cal B$ be the basin of attraction and 
$\Cal O$ be the operating range. It is more reasonable to measure 
the size of the subset $\Cal N\subset\Cal  B$ consisting of all states 
whose forward trajectories remain in the operating range, see Figure 3.
We will call $\Cal N$ {\it the normal operating range}. 
It is just as difficult to form a reasonable 
measure of the size of $\Cal  N$.

\bigskip


\hskip 100bp\psfig{file=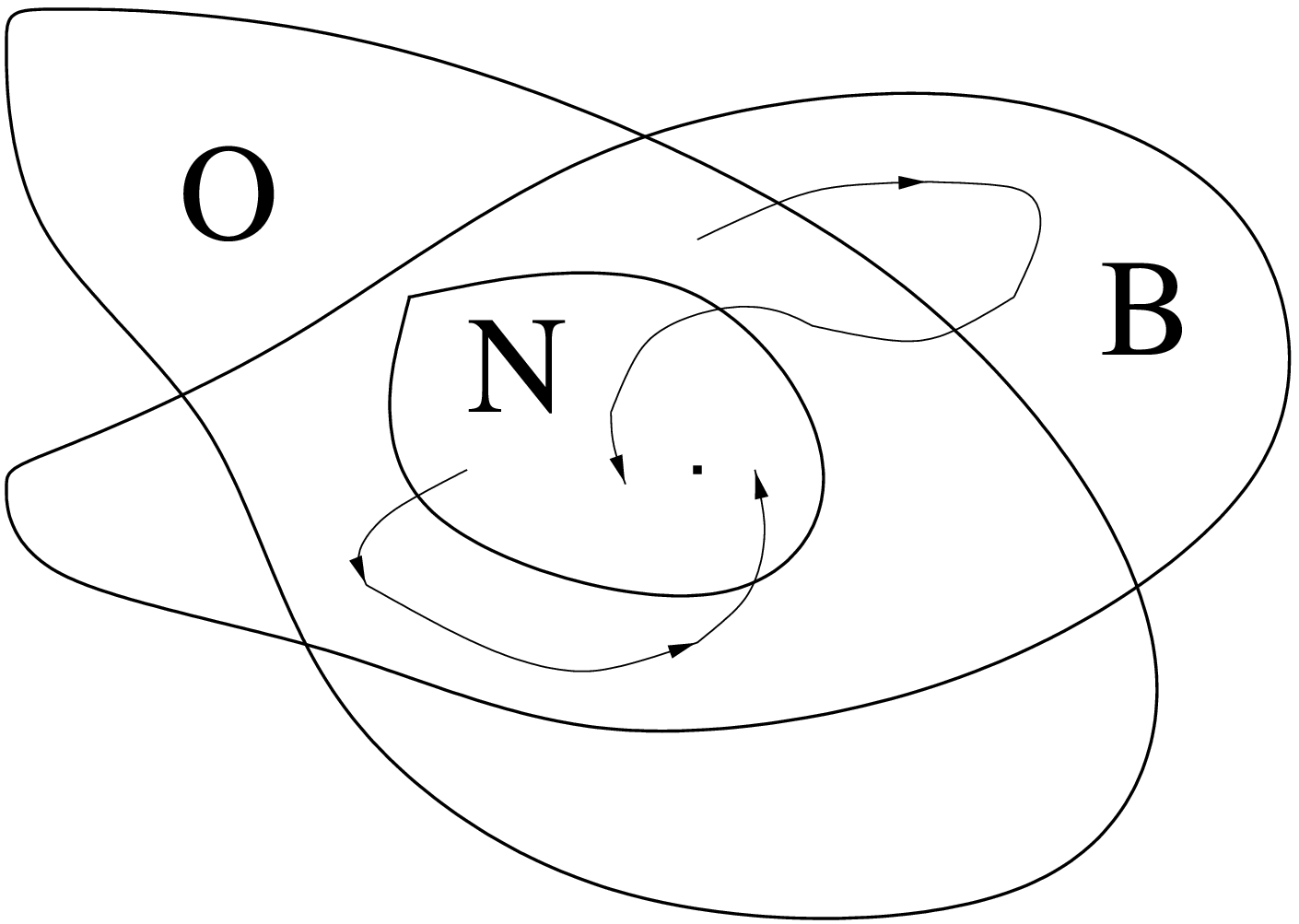,width=3.5truein,height=2truein}%

\centerline{Figure 3}
\smallskip

\proclaim{Problem} Find an effective analytical or numerical method to 
estimate the injectivity radius of 
the normal operating range or define a better 
measure of ``size" together with an
effective analytical or numerical estimate.
\endproclaim

In engineering applications the time it takes 
to drive the states to the equilibrium is very important. 
This is usually characterized by the rise time, peak time, 
settling time, etc., \cite{11, p. 222}. 
These notions are only well defined 
for solutions of
$\;
m\,\ddot z + c\,\dot z +k\,z = f\, 
$
with $m$, $c$, $k$, and $f$ constants. 
Considering the initial conditions with $\dot z=0$ only, 
one rescales the initial value problem to 
$\;
\ddot z + 2\zeta\,\dot z + z = 1$ ,$\;z(0)=0,\;\dot z(0)=0$. 
The rise time, $T_r$, is the first time $z(t)$ reaches $1$. 
The time constant is $1/\zeta$. The settling time is the time 
it takes $z(t)$ to get and stay within 2\% of the final value $1$. 
These notions are used for linear constant coefficient ODEs 
which have two dominant poles. 

One can define analogues of these notions for nonlinear systems.
An engineer may choose a {\it target operating range\/} 
$\Cal D$, which is a small neighborhood of the equilibrium. 

\bigskip


\hskip 100bp\psfig{file=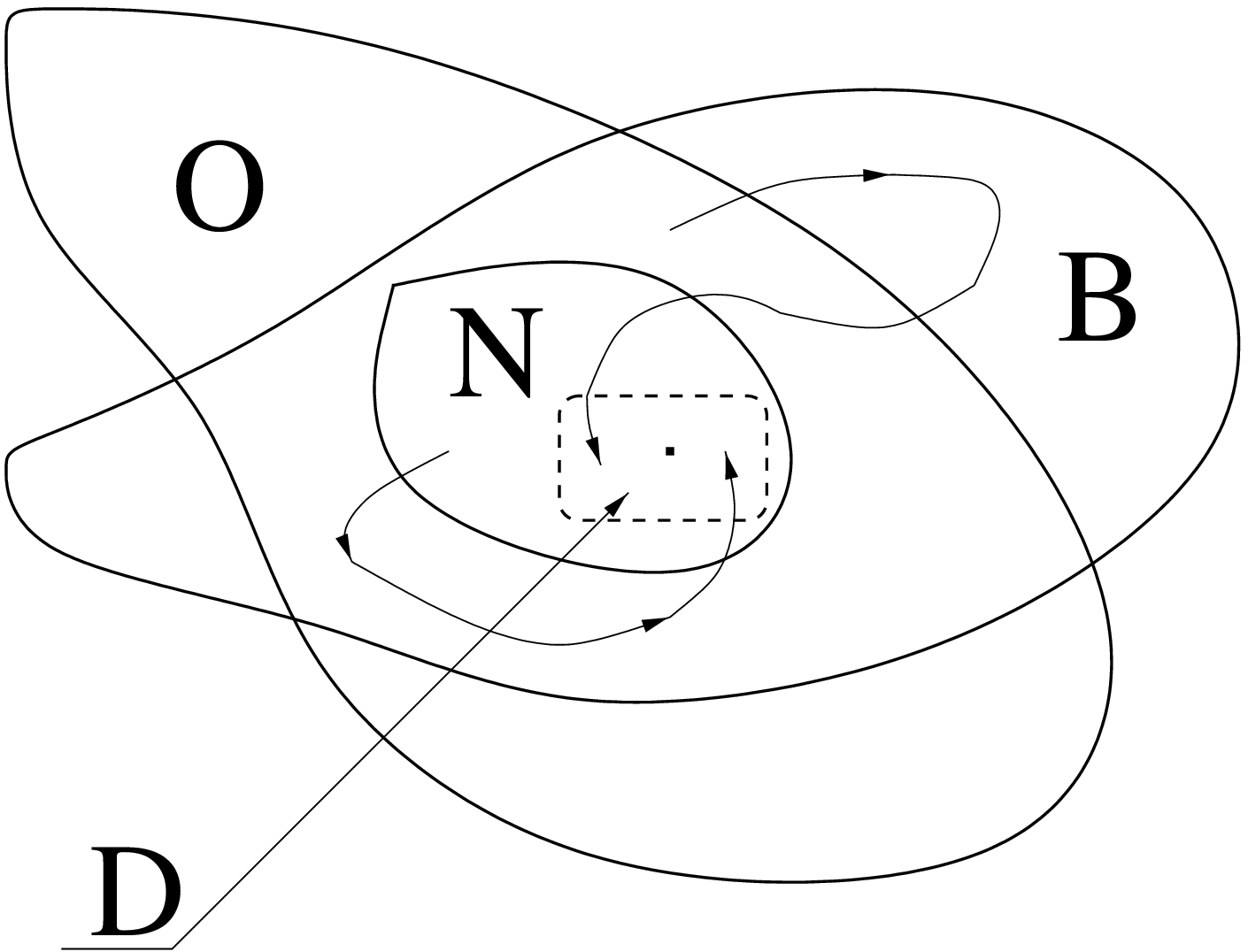,width=3.5truein,height=2truein}%

\centerline{Figure 4}
\smallskip

\noindent Define 
the $\Cal N\Cal D$-settling time 
$$
T_{\Cal N\Cal D}\,=\,\sup_{x_0\in\Cal N} 
\{t\ge 0 \vert \; x(t;x_0)\notin \Cal D\}\,.
$$
This is the time after which any trajectory starting in 
the normal operating range 
$\Cal N$ will settle into the target range $\Cal D$. 

\proclaim{Problem} Find an effective analytical or numerical method to 
estimate the $\Cal N\Cal D$-settling time or define a better 
measure of ``short time" together with 
an effective analytical or numerical estimate.
\endproclaim

Another important characteristic of a control law is the cost. 
One often wishes to minimize some function of the trajectory 
and/or the control input (for example, minimize the energy used 
to complete a task). This imposes a non-trivial 
restriction on a control problem.
Real life systems have additional physical restrictions. For example, 
there is a maximal voltage which may be applied to a DC motor before 
it saturates. 

Reiterating, a good stabilizing control law should maximize the basin of 
attraction and minimize the time and cost while operating within 
the physical restrictions of the given system.

\bigskip

The problem of maximizing the size of the basin of attraction 
taken to the extreme leads to the question of finding the 
topologically best control laws. For example, the state space 
for the inverted pendulum cart is $S^1\times\Bbb R\times\Bbb R^2$. 
For topological reasons it is impossible to find a 
control law so that the resulting flow has a globally 
asymptotically stable equilibrium in this case. 
The best one could hope for is a flow with a regular 
compact global attractor, \cite{12, 16}. Recall that such an attractor 
is the union of the unstable manifolds of finitely many 
hyperbolic fixed points. Also recall that the shape 
(in the sense of K. Borsuk) of the 
(global) attractor must be the same as the shape of the state space, 
\cite{15}.
We wish to minimize the number of 
fixed points. For the inverted pendulum cart 
{\it the topologically best flow} 
is depicted in Figure 5.

\bigskip 


\hskip 90bp\psfig{file=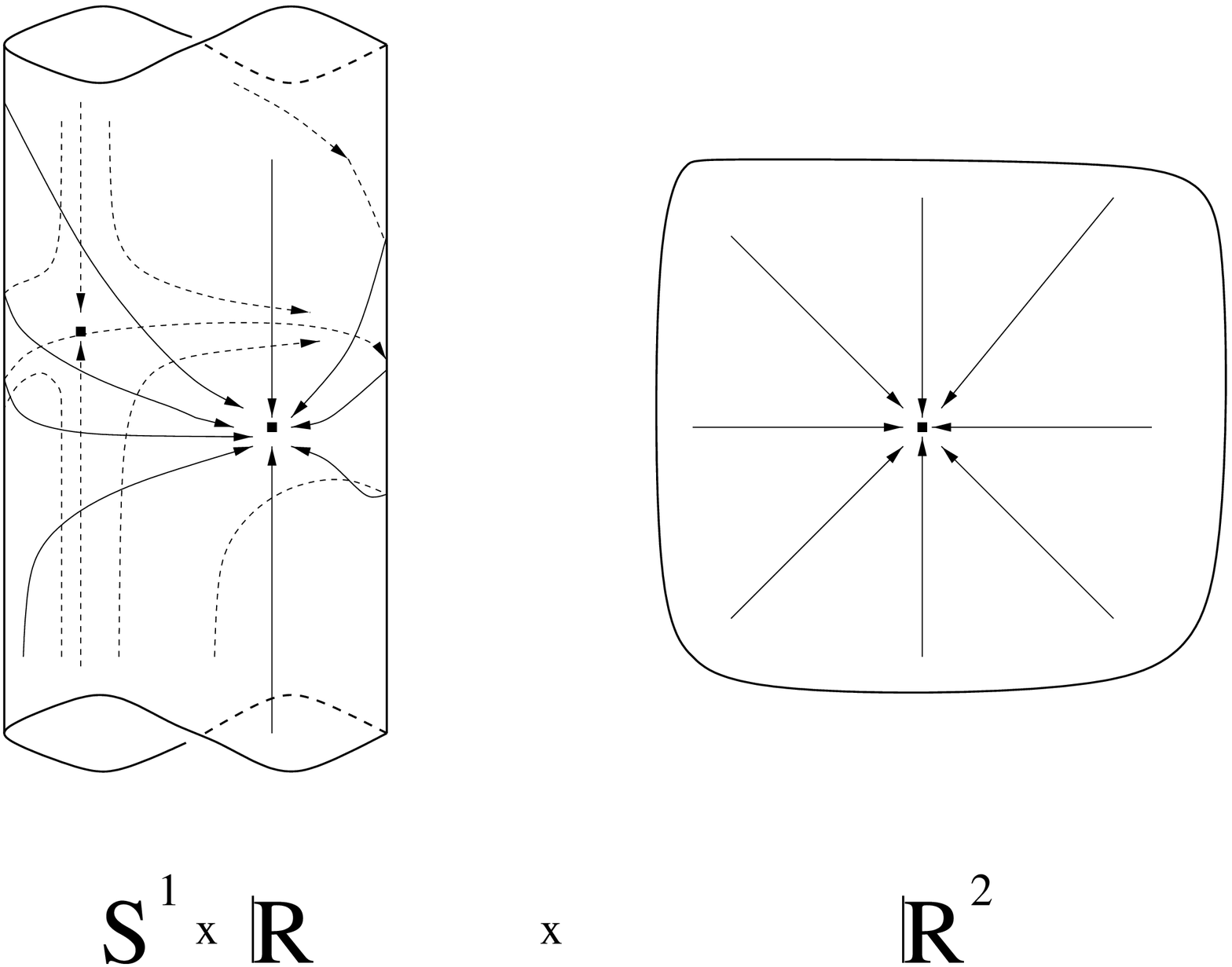,width=4truein,height=2.6truein}%

\centerline{Figure 5}
\smallskip

\noindent This flow must have one index 0 fixed point 
and one index 1 fixed point. 
The basin of attraction of the index 0 fixed point is the complement 
of a properly embedded $\Bbb R^3$.

\proclaim{Problem} When is it possible to find a control input 
producing a topologically best flow? Is there a method to find such 
control laws?
\endproclaim

At the moment we do not even know whether a topologically best 
control law exists for the inverted pendulum cart. Furthermore, 
we do not know whether a control law with a global attractor 
exists!
\bigskip

A related problem is to describe the behavior of the system 
far away from the equilibrium. It is possible to compactify 
some dynamical systems with algebraic nonlinearities, \cite{10}. 
The advantage of a compactification is that overall dynamics 
may be well understood by linearization about each of 
the critical submanifolds. The problem is that for 
nonalgebraic nonlinearities it is not known when a compactification 
exists. It is even interesting and instructive to try to 
compactify the system $\ddot\theta+c\,\dot\theta+\sin(\theta)=0$.
\bigskip

We wrote this paper to encourage other people to think about 
nonlinear control theory.  
It has a great many unanswered fundamental questions. 
The field is ready for new ideas. 
\bigskip

\hyphenation{Sprin-g-er-Ver-lag}

\widestnumber\key{442}

\centerline {REFERENCES}
\medskip

\ref\no 1\by D. Auckly, L. Kapitanski, A. Kelkar, W. White
\paper Matching and pole placement for underactuated systems
\paperinfo Preprint
\yr 1999
\endref

\ref\no 2\by D. Auckly, L. Kapitanski, A. Kelkar, W. White
\paper Matching and digital control implementation for 
underactuated systems
\paperinfo Preprint
\yr 1999
\endref

\ref\no 3\by D. Auckly, L. Kapitanski, W. White
\paper Control of nonlinear underactuated systems
\jour to appear in Communications on Pure Appl. Math.
\yr 1999
\endref

\ref\no 4\by N. P. Bhatia and G. P. Szeg\"o
\book Dynamical systems: stability theory
and applications
\bookinfo Lecture Notes in Mathematics, 35
\publ Springer-Verlag
\publaddr Berlin Heidelberg New York
\yr 1967
\endref

\ref\no 5\by A. Bloch, N. Leonard and J. Marsden
\paper Stabilization of mechanical systems using controlled 
Lagrangians
\inbook Proc. IEEE Conf. Dec. Contr.
\publ San Diego
\publaddr Calif.
\yr 1997
\pages 2356-2361
\endref

\ref\no 6\by A. Bloch, N. Leonard and J. Marsden
\paper Matching and stabilization by the method of controlled 
Lagrangians
\inbook Proc. IEEE Conf. Dec. Contr.
\yr 1998
\publ Tampa
\publaddr Fla.
\pages 1446-1451
\endref

\ref\no 7\by A. Bloch, N. Leonard and J. Marsden
\paper Stabilization of the pendulum on a rotor arm 
by the method of controlled Lagrangians
\inbook  Proc. IEEE Int. Conf. on Robotics and Automation
\yr 1999
\publ Detroit
\publaddr Mich.
\pages 500-505
\endref

\ref\no 8\by A. Bloch, N. Leonard and J. Marsden
\paper Controlled Lagrangians and a stabilization of mechanical systems 
 I: The first matching theorem
\jour Preprint
\yr 1999
\endref

\ref\no 9\by A. Bloch, N. Leonard and J. Marsden
\paper Potential shaping and the method of controlled Lagrangians
\jour Preprint
\yr 1999
\endref

\ref\no 10\by O. I. Bogoyavlensky
\book Methods in the qualitative theory of dynamical systems 
in astrophysics and gas dynamics
\bookinfo Springer Series in Soviet Mathematics 
\publ Springer-Verlag
\publaddr Berlin-New York
\yr 1985
\endref

\ref\no 11\by R. C. Dorf and R. H. Bishop
\book Modern control systems
\bookinfo
\publ Addison-Wesley
\yr 1995
\endref

\ref\no 12\by J. K. Hale
\book Asymptotic Behavior of Dissipative Systems
\publ Amer. Math. Soc.
\publaddr Providence, RI
\bookinfo Mathematical Surveys and Monographs, vol. 25
\yr 1987
\endref

\ref\no 13\by J. Hamberg
\paper General matching conditions in the theory of 
controlled Lagrangians
\inbook to appear in proceedings of 
the 38th Conference on Decision and Control
\publ Phoenix
\publaddr Ariz.
\yr 1999
\endref

\ref\no 14\by T. Kailath
\book Linear systems
\publ Prentice-Hall
\publaddr Englewood Cliffs, N.J.
\yr 1980
\endref

\ref\no 15\by L. Kapitanski and I. Rodnianski
\paper Shape and Morse theory of attractors
\jour Communications in Pure Appl. Math.
\yr 1999
\toappear
\endref

\ref\no 16\by O. A. Ladyzhenskaya
\book Attractors for Semigroups and
Evolution Equations
\bookinfo (Lezioni Lin\-cee)
\publ Cambridge University Press
\publaddr Cambridge
\yr 1991
\endref

\ref\no 17\by A. J. van der Schaft
\paper Stabilization of Hamiltonian systems
\jour Nonlinear Analysis, Theory, Methods \& Applications
\vol 10 \issue 10
\pages 1021-1035
\yr 1986
\endref

\ref\no 18\by E. D. Sontag
\book Mathematical control theory
\bookinfo Texts in Applied Mathematics, 6
\publ Springer-Verlag
\publaddr New York
\yr 1990
\endref

\end